\documentclass[11pt,twoside]{article}

\usepackage{amsfonts,epsfig,graphicx,subfigure}
\usepackage{amsmath,amssymb,amsthm}

\usepackage{epic,eepic}

\usepackage{epsf}
\usepackage{fancyheadings}
\usepackage{graphics}
\usepackage{amsfonts}
\usepackage{amsmath}
\usepackage{amssymb}
\usepackage{epic}
\usepackage{eepic}
\usepackage{eepicemu}

\newenvironment{carlist}
 {\begin{list}{$\bullet$}
 {\setlength{\topsep}{0in} \setlength{\partopsep}{0in}
  \setlength{\parsep}{0in} \setlength{\itemsep}{\parskip}
  \setlength{\leftmargin}{0.07in} \setlength{\rightmargin}{0.08in}
  \setlength{\listparindent}{0in} \setlength{\labelwidth}{0.08in}
  \setlength{\labelsep}{0.1in} \setlength{\itemindent}{0in}}}
 {\end{list}}

\newcommand{\bcar}{\begin{carlist}}
\newcommand{\ecar}{\end{carlist}}

\newenvironment{carliste}
 {\begin{list}x
 {\setlength{\topsep}{0in} \setlength{\partopsep}{0in}
  \setlength{\parsep}{0in} \setlength{\itemsep}{\parskip}
  \setlength{\leftmargin}{0.07in} \setlength{\rightmargin}{0.08in}
  \setlength{\listparindent}{0in} \setlength{\labelwidth}{0.08in}
  \setlength{\labelsep}{0.1in} \setlength{\itemindent}{0in}}}
 {\end{list}}

\newcommand{\bcare}{\begin{carliste}}
\newcommand{\ecare}{\end{carliste}}

\newcommand{\Prob}{\ensuremath{\mathbb{P}}}

\makeatletter
\long\def\@makecaption#1#2{
        \vskip 0.8ex
        \setbox\@tempboxa\hbox{\small {\bf #1:} #2}
        \parindent 1.5em  
        \dimen0=\hsize
        \advance\dimen0 by -3em
        \ifdim \wd\@tempboxa >\dimen0
                \hbox to \hsize{
                        \parindent 0em
                        \hfil 
                        \parbox{\dimen0}{\def\baselinestretch{0.96}\small
                                {\bf #1.} #2
                                } 
                        \hfil}
        \else \hbox to \hsize{\hfil \box\@tempboxa \hfil}
        \fi
        }
\makeatother

\long\def\comment#1{}

\makeatletter
\def\@cite#1#2{[\if@tempswa #2 \fi #1]}
\makeatother

\makeatletter
\long\def\barenote#1{
    \insert\footins{\footnotesize
    \interlinepenalty\interfootnotelinepenalty 
    \splittopskip\footnotesep
    \splitmaxdepth \dp\strutbox \floatingpenalty \@MM
    \hsize\columnwidth \@parboxrestore
    {\rule{\z@}{\footnotesep}\ignorespaces
      #1\strut}}}
\makeatother

\newcommand{\Ind}{\ensuremath{\mathbb{I\,}}}

\newcommand{\bit}{\begin{itemize}}
\newcommand{\eit}{\end{itemize}}
\newcommand{\ben}{\begin{enumerate}}
\newcommand{\een}{\end{enumerate}}

\newcommand{\bear}{\begin{eqnarray}}
\newcommand{\eear}{\end{eqnarray}}

\newcommand{\vtiny}{\vspace*{.1in}}

\newcommand{\df}{\ensuremath{d}}

\newcommand{\prob}{\ensuremath{{\mathbb{P}}}}

\newcommand{\Exs}{\ensuremath{{\mathbb{E}}}}

\newcommand{\beq}{\begin{quotation}}
\newcommand{\enq}{\end{quotation}}

\newcommand{\estart}{\begin{equation}}
\newcommand{\eend}{\end{equation}}

\newcommand{\edist}{\ensuremath{\overset{d}{=}}}

\newcommand{\defn}{\ensuremath{:  =}}

\newcommand{\Ysca}{{{Y}}}

\newcommand{\Wsca}{{{W}}}

\newcommand{\bec}{\begin{center}}
\newcommand{\enc}{\end{center}}

\newcommand{\beit}{\begin{itemize}}
\newcommand{\enit}{\end{itemize}}

\newcommand{\been}{\begin{enumerate}}
\newcommand{\enen}{\end{enumerate}}

\newcommand{\comsl}{\begin{slide}}
\newcommand{\comspor}{\begin{slide*}}

\newcommand{\comsld}[2]{\begin{slide}[#1,#2]}
\newcommand{\comspord}[2]{\begin{slide*}[#1,#2]}

\newcommand{\mendsl}{\end{slide}}
\newcommand{\mendspo}{\end{slide*}}

\newcommand{\estim}[1]{\ensuremath{\widehat{#1}}}

\newcommand{\real}{\ensuremath{{\mathbb{R}}}}

\theoremstyle{plain}

\newtheorem{theo}{Theorem}[section]

\newtheorem{lem}{Lemma}[section]
\newtheorem{prop}{Proposition}[section]
\newtheorem{cor}{Corollary}[section]

\theoremstyle{definition} 

\newtheorem{nota}{Notation}[section]
\newtheorem{de}{Definition}[section]
\newtheorem{exa}{Example}[section]
\newtheorem{as}{Assumption}[section]
\newtheorem{alg}{Algorithm}[section]

\newcommand{\btheo}{\begin{theo}}
\newcommand{\bde}{\begin{de}}
\newcommand{\ble}{\begin{lem}}
\newcommand{\bpr}{\begin{prop}}
\newcommand{\bno}{\begin{nota}}
\newcommand{\bex}{\begin{exa}}
\newcommand{\bcor}{\begin{cor}}
\newcommand{\spro}{\begin{proof}}
\newcommand{\bas}{\begin{as}}
\newcommand{\balg}{\begin{alg}}

\newcommand{\etheo}{\end{theo}}
\newcommand{\ede}{\end{de}}
\newcommand{\ele}{\end{lem}}
\newcommand{\epr}{\end{prop}}
\newcommand{\eno}{\end{nota}}
\newcommand{\eex}{\end{exa}}
\newcommand{\ecor}{\end{cor}}
\newcommand{\fpro}{\end{proof}}
\newcommand{\eas}{\end{as}}
\newcommand{\ealg}{\end{alg}}

\theoremstyle{plain}

\newtheorem{theos}{Theorem}
\newtheorem{props}{Proposition}
\newtheorem{lems}{Lemma}
\newtheorem{cors}{Corollary}

\theoremstyle{definition}
\newtheorem{exas}{Example}
\newtheorem{algs}{Algorithm}
\newtheorem{asss}{Asumption}
\newtheorem{defns}{Definition}

\newcommand{\btheos}{\begin{theos}}
\newcommand{\etheos}{\end{theos}}
\newcommand{\bprops}{\begin{props}}
\newcommand{\eprops}{\end{props}}
\newcommand{\bdes}{\begin{defns}}
\newcommand{\edes}{\end{defns}}
\newcommand{\blems}{\begin{lems}}
\newcommand{\elems}{\end{lems}}
\newcommand{\bcors}{\begin{cors}}
\newcommand{\ecors}{\end{cors}}
\newcommand{\bexs}{\begin{exas}}
\newcommand{\eexs}{\end{exas}}
\newcommand{\balgs}{\begin{algs}}
\newcommand{\ealgs}{\end{algs}}
\newcommand{\bass}{\begin{asss}}
\newcommand{\eass}{\end{asss}}

\newcommand{\Sset}{\ensuremath{S}}

\newcommand{\Amat}{\ensuremath{X}}
\newcommand{\Amatt}[1]{\ensuremath{\Amat_{#1}}}
\newcommand{\arow}{\ensuremath{x}}
\newcommand{\acol}{\ensuremath{X}}

\newcommand{\betastar}{\ensuremath{\beta^*}}

\newcommand{\Esca}{\ensuremath{W}}

\newcommand{\numobs}{\ensuremath{n}}
\newcommand{\spindex}{\ensuremath{s}}
\newcommand{\mdim}{\ensuremath{p}}

\newcommand{\Ker}{\ensuremath{\operatorname{Ker}}}
\newcommand{\Range}{\ensuremath{\operatorname{Ra}}}

\newcommand{\Uset}{\ensuremath{U}}
\newcommand{\Vset}{\ensuremath{V}}

\newcommand{\Beta}[1]{\ensuremath{\beta_{#1}}}

\newcommand{\OptFun}{\ensuremath{f}}

\newcommand{\Sbest}{\ensuremath{\estim{\Sset}}}

\newcommand{\MatInv}[1]{\ensuremath{\left[ #1^T #1 \right]^{-1}}}

\newcommand{\Gfullmat}[1]{\ensuremath{#1 \MatInv{#1} #1^T}}

\newcommand{\Metric}{\ensuremath{\Delta}}

\newcommand{\SbackU}{\ensuremath{S \backslash U}}
\newcommand{\UbackS}{\ensuremath{U \backslash S }}

\newcommand{\pdim}{\mdim}

\newcommand{\myones}{\vec{1}}

\newcommand{\mybetavec}{\ensuremath{\vec{v}}}

\newcommand{\noncent}{\ensuremath{\nu}}

\newcommand{\Numset}{\ensuremath{N}}

\newcommand{\perr}{\ensuremath{p_{\operatorname{err}}}}

\newcommand{\mprob}{\ensuremath{\mathbb{P}}}

\newcommand{\decopt}{\ensuremath{\phi_{\operatorname{opt}}}}

\newcommand{\goodW}{\ensuremath{\mathcal{A}}}

\newcommand{\Projmat}[1]{\ensuremath{\Pi_{#1}}}
\newcommand{\Projmatc}[1]{\ensuremath{\Pi^\perp_{#1}}}

\newcommand{\minval}{\ensuremath{\mathcal{M}}}

\newcommand{\mycfrac}{\ensuremath{\gamma}}

\newcommand{\gamvar}{\ensuremath{\gamma}}

\newcommand{\myparagraph}[1]{\paragraph{#1:}}

\begin{document}
\begin{center}

{{\LARGE \bf{Information-theoretic limits on sparsity recovery in the
high-dimensional and noisy setting}}}

\vspace*{.5in}

{\large {
\begin{center}
Martin J. Wainwright \\
Department of Statistics, and \\
Department of Electrical Engineering and Computer Sciences \\
University of California, Berkeley \\
\texttt{wainwrig@\{eecs,stat\}.berkeley.edu}
\end{center}

}}

\vtiny

{\large Technical Report, UC Berkeley, Department of Statistics} \\
{\large January 2007}

\vtiny

\begin{abstract}
The problem of recovering the sparsity pattern of a fixed but unknown
vector $\betastar \in \real^\pdim$ based on a set of $\numobs$ noisy
observations arises in a variety of settings, including subset
selection in regression, graphical model selection, signal denoising,
compressive sensing, and constructive approximation.  Of interest are
conditions on the model dimension $\pdim$, the sparsity index
$\spindex$ (number of non-zero entries in $\betastar$), and the number
of observations $\numobs$ that are necessary and/or sufficient to
ensure asymptotically perfect recovery of the sparsity pattern.  This
paper focuses on the information-theoretic limits of sparsity
recovery: in particular, for a noisy linear observation model based on
measurement vectors drawn from the standard Gaussian ensemble, we
derive both a set of sufficient conditions for asymptotically perfect
recovery using the optimal decoder, as well as a set of necessary
conditions that \emph{any} decoder, regardless of its computational
complexity, must satisfy for perfect recovery.  This analysis of
optimal decoding limits complements our previous
work~\cite{Wainwright06a_aller} on sharp thresholds for sparsity
recovery using the Lasso ($\ell_1$-constrained quadratic programming)
with Gaussian measurement ensembles.
\end{abstract}
\end{center}



\section{Introduction}

Suppose that we are given a set of $\numobs$ observations of a fixed
but unknown vector $\betastar \in \real^\pdim$.  In a variety of
settings, it is known \emph{a priori} that the vector $\betastar$ is
sparse, meaning that its support set $\Sset$---corresponding to those
indices $i$ for which $\betastar_i$ is non-zero---is relatively small,
say with \mbox{size $|\Sset| =: \, \spindex \ll \pdim$}.  Sparsity
recovery refers to the problem of correctly estimating the support set
$\Sset$ based on a set of noisy observations.  This sparsity recovery
problem is of broad interest, arising in various areas, including
subset selection in regression~\cite{Miller90}, structure estimation
in graphical models~\cite{Meinshausen06}, sparse
approximation~\cite{Devore93,Natarajan95}, signal
denoising~\cite{Chen98}, and compressive
sensing~\cite{Donoho06,CandesTao05}.

A great deal of work over the past few years has focused on the
performance of computationally tractable methods, many based on
$\ell_1$ or other convex relaxations, both for recovering the exact
sparsity pattern as well as related problems in sparse approximation.
We provide a brief overview of those parts of this extensive
literature most relevant to our work in Section~\ref{SecRelatedWork}
below.  Of equal interest and complementary in nature, however, are
the information-theoretic limits associated with the performance of
\emph{any} procedure for sparsity recovery.  Such understanding of
fundamental limitations is crucial in assessing the behavior of
computationally tractable methods.  In particular, there is little
point in proposing novel methods for sparsity recovery, possibly with
higher computational complexity, if currently extant and
computationally tractable methods achieve the information-theoretic
limits.  On the other hand, an information-theoretic analysis can
reveal where there currently exists a gap between the performance of
computationally tractable methods, and the fundamental limits. Indeed,
the information-theoretic analysis of this paper makes contributions
of both types.

With this motivation in mind, the focus of this paper is on the
information-theoretic limitations of sparsity recovery.  In
particular, our analysis focuses on the noisy and high-dimensional
setting, meaning that the observations are contaminated by noise, and
all three problem parameters---the \emph{number of observations}
$\numobs$, the \emph{model dimension} $\pdim$, and the \emph{sparsity
index} $\spindex$, defined below---may tend to infinity.  Our main
results, stated more precisely in Section~\ref{SecMainResults}, are
necessary and sufficient conditions on the triplet $(\numobs, \pdim,
\spindex)$ for exact recovery.  In particular, given noisy linear
observations based on measurement vectors drawn from the standard
Gaussian ensemble, we derive both a set of sufficient conditions for
asymptotically perfect recovery using the optimal decoder, as well as
a set of necessary conditions that any decoder must satisfy for
perfect recovery.  The analysis given here complements our earlier
paper~\cite{Wainwright06a_aller} that established precise thresholds
on the success/failure of the Lasso (i.e., $\ell_1$-constrained
quadratic programming) for sparsity recovery.

The remainder of this paper is organized as follows.  In
Section~\ref{SecRelatedWork}, we provide a more precise formulation of
the problem, and a brief discussion of past work, whereas
Section~\ref{SecMainResults} provides a precise statement of our main
results, and a discussion of their consequences.
Section~\ref{SecAnalysis} and the appendices are devoted to the proofs
of our main results, and we conclude in Section~\ref{SecDiscussion}
with a discussion of open directions.

\subsection{Problem formulation and past work}
\label{SecRelatedWork}

We begin with a more precise formulation of the problem, as well as a
discussion of previous work, with emphasis on that most closely
related to the results in this paper. Let $\betastar \in \real^\pdim$
be a fixed but unknown vector; we refer to the ambient dimension
$\pdim$ as the \emph{model dimension}.  Define the support set of
$\betastar$ as
\begin{eqnarray}
\Sset & \defn & \{ i \in \{1, \ldots, \mdim \} \; \mid \; \betastar_i
\neq 0 \}.
\end{eqnarray}
We refer to its size $\spindex \defn |\Sset|$ as the \emph{sparsity
index}.  Finally, suppose that we are given a set of $\numobs$
observations, of the form
\begin{eqnarray}
\label{EqnLinearObs}
\Ysca_i & = & \arow_i^T \betastar + \Wsca_i,  \qquad i = 1, \ldots, \numobs
\end{eqnarray}
where each $\arow_i \in \real^\pdim$ is a measurement vector, and
$\Wsca_i \sim N(0, \sigma^2)$ is additive Gaussian noise.  Of interest
are conditions on the triplet $(\numobs, \pdim, \spindex)$ under which
a given method either succeeds or fails in recovering the sparsity
pattern $\Sset$.  

\myparagraph{Observation models} The linear observation
model~\eqref{EqnLinearObs} can be studied in either its noiseless
variant ($\sigma^2 = 0$), or the noisy setting ($\sigma^2 > 0$); this
paper focuses exclusively the noisy setting.  In addition, previous
work has addressed both deterministic families and random ensembles of
measurement vectors $\{\arow_i \}_{i=1}^\numobs$.  The analysis in
this paper is based on the \emph{standard Gaussian measurement
ensemble}, in which each measurement vector $\arow_i$ is drawn from
the zero-mean isotropic Gaussian distribution $N(0, I_{\pdim \times
\pdim})$.

\myparagraph{Error metrics} Consider some method that generates the
vector $\estim{\beta} \in \real^\pdim$ as an estimate of the truth
$\betastar$.  There are various distinct criteria for assessing how
close the estimate is to the truth, including
\bit
\item various $\ell_p$ norms $\Exs \|\estim{\beta} - \betastar \|_p$,
especially $\ell_2$ and $\ell_1$, or
\item some measurement of predictive power (e.g., $\Exs[\|Y_i -
\estim{Y}_i \|_2^2]$, where $\estim{Y}_i$ is the estimate based on
$\estim{\beta}$).
\eit 
Given the abundance of recent results on sparse approximation (not all
of which are mutually comparable), it is particularly important to
specify up front the choice of error metric.  In this paper, we focus
exclusively on the sparsity recovery problem, for which the
appropriate error metric is simply the $0-1$ loss associated with the
event of recovering the correct support $\Sset$---viz.:
\begin{eqnarray}
\label{EqnRecover}
\rho(\estim{\beta}, \betastar) & = & \Ind \left[ \left \{
\estim{\beta}_i \neq 0 \quad \forall i \in \Sset \right \} \cap \left
\{ \estim{\beta}_j = 0 \quad \forall j \notin \Sset \right \} \right].
\end{eqnarray}

\myparagraph{Past work} Closely related in its information-theoretic
spirit is the earlier paper of Fletcher et al.~\cite{Fletcher06} that
analyzed the standard Gaussian ensemble from a rate-distortion
perspective, studying the average $\ell_2$-error of the optimal
decoder.  The results given here also address the
information-theoretic limitations, albeit of the sparsity recovery
problem, using the error metric~\eqref{EqnRecover} as opposed to
$\ell_2$-norm.  In a related but distinct line of work, the use of
$\ell_1$-relaxation for sparse approximation has a lengthy history;
relatively early papers from the 1990s include the work of Chen,
Donoho and Saunders~\cite{Chen98}, as well as
Tibshirani~\cite{Tibshirani96} on $\ell_1$-constrained quadratic
programming (known as the Lasso in the statistics literature).  A
great deal of subsequent work has analyzed the performance of
$\ell_1$-relaxations, both in the
noiseless~\cite{Elad02,Feuer03,Malioutov04} and noisy
setting~\cite{Tropp06} for deterministic ensembles, as well as the
noiseless~\cite{Donoho04a,CandesTao05,DonTan06} and noisy
setting~\cite{CandesTao06,CanRomTao06,Donoho04b,Meinshausen06,Zhao06,Wainwright06a_aller}
for random ensembles.  Other work has provided conditions under which
estimation of a noise-contaminated vector via the
Lasso~\cite{CanRomTao06,Donoho04b} or other types of convex
relaxation~\cite{CandesTao06} is stable in the $\ell_2$ sense;
however, such $\ell_2$-stability does not guarantee exact recovery of
the underlying sparsity pattern.

A notable feature of the results given here is that they apply to
completely general scaling of the triplet $(\numobs, \pdim,
\spindex)$.  In contrast, most previous work has addressed one of two
possible special cases of sparsity scaling: (a) either the
\emph{linear sparsity regime}~\cite[e.g]{
CandesTao05,Donoho04a,Donoho04b}, in which $\spindex = \alpha \pdim$
for some $\alpha \in (0,1)$; or (b) the \emph{sublinear sparsity
regime}~\cite[e.g.,]{Meinshausen06,Zhao06}, in which $\spindex/\pdim$
tends to zero.  Depending on the underlying motivation for sparse
approximation, both of these sparsity regimes are of independent
interest. In covering the full range of scaling, the results given
here are complementary to those of our previous
paper~\cite{Wainwright06a_aller} that provided threshold results, also
applicable to general scaling of $(\numobs, \pdim, \spindex)$, for the
success/failure of the Lasso when used for sparsity recovery with
random Gaussian measurement ensembles.  We discuss connections to
previous work in more technical detail following the statement of our
main results below.

\subsection{Our contributions}
\label{SecMainResults}

The analysis of this paper procedure is asymptotic in nature, focusing
on scaling conditions on the triplet $(\numobs, \mdim, \spindex)$
under which asymptotically exact recovery is either possible or
impossible.  As mentioned previously, we focus on the linear
observation model~\eqref{EqnLinearObs} in the noisy setting ($\sigma^2
> 0$), and with the measurement vectors $\arow_i$ drawn in an i.i.d.
manner from the standard Gaussian $N(0, I_{\pdim \times \pdim})$
ensemble. A decoder is a mapping from the $\numobs$-vector of
observations $\Ysca$ to an estimated subset---say of the form $\Sbest
= \phi(Y)$.  We think of the underlying true vector $\betastar \in
\real^\pdim$ with its support $\Sset$ randomly chosen, uniformly over
all ${\pdim \choose \spindex}$ subspaces of size $\spindex$.
Accordingly, the average error probability $\perr$ of any decoder is
given by
\begin{equation*}
\perr(\phi) \; = \; \frac{1}{{\pdim \choose \spindex}} \sum_{\Sset, \;
|\Sset| = \spindex} \Prob[\phi(Y) \neq \Sset \, \mid \, \Sset].
\end{equation*}
Here the term $\mprob[\phi(Y) \neq \Sset \, \mid \, \Sset]$
corresponds to the probability, conditioned on the true underlying
support being $\Sset$ and averaging over the measurement noise
$\Wsca$, the choice of Gaussian random matrix $\Amat$, and the choice
of the entries $\beta^*_{\Sset}$ on the fixed support $\Sset$, that
the decoder makes an error.  We say that
\begin{itemize}
\item the sparsity recovery is \emph{asymptotically reliable}
(error-free) if $\perr(\phi) \rightarrow 0$ as $\numobs \rightarrow
+\infty$, and
\item the sparsity recovery is \emph{asymptotically unreliable} if for
  some constant $c > 0$, the error probability stays bounded
  $\perr(\phi) \geq c$ as $\numobs \rightarrow +\infty$.
\end{itemize}
\noindent In addition to the three parameters $(\numobs, \mdim,
\spindex)$, our results also involve the minimum value of the unknown
vector $\betastar$ on its support, given by
\begin{eqnarray}
\label{EqnDefnMinVal}
\minval(\betastar) & \defn & \min_{i \in \Sset} |\beta^*_i|.
\end{eqnarray}
We begin by stating a set of conditions on the triplet $(\numobs,
\pdim, \spindex)$ which are sufficient to ensure asymptotically
perfect recovery of the sparsity pattern:
\btheos[Sufficient conditions]
\label{ThmSuff}
If $(\numobs - \spindex) \minval^2(\betastar) \rightarrow +\infty$,
then the following condition suffices to ensure asymptotically
reliable recovery: for some fixed constant $C > 0$,
\begin{eqnarray}
\label{EqnSuffGeneral}
\numobs & > & C \; \max\left \{ \spindex \log (\pdim/\spindex), \;
\frac{1}{\minval^2(\betastar)} \log(\pdim - \spindex) \right \}.
\end{eqnarray}
\etheos
\noindent The proof of this claim, given in Section~\ref{SecSuff}, is
constructive in nature, based on direct analysis of the error
probability associated with the optimal decoder.
\btheos[Necessary conditions]
\label{ThmNec}
Asymptotically reliable recovery is impossible under the following
condition: for some fixed constant $C' > 0$:
\begin{eqnarray}
\label{EqnNecGeneral}
\numobs & < & \left[\frac{C'}{\spindex \; \minval^2(\betastar)}
\right] \; \spindex \log \frac{\pdim}{\spindex}.
\end{eqnarray}
\etheos
\noindent The proof of this claim, given in Section~\ref{SecNec}, is
somewhat more indirect in nature, based on exploiting a corollary of
Fano's inequality~\cite{Cover, Hasminskii78,IbrHas81,Yu97}, in order
to lower bound the probability of error for a restricted hypothesis
testing problem.  To interpret these results, we consider two distinct
regimes of sparsity:
\myparagraph{Regime of sublinear sparsity} First suppose that the
sparsity is sublinear, meaning that \mbox{$\spindex = o(\pdim)$.}
Based on the two theorems, we identify the critical scaling as
$\minval^2(\betastar) = \Theta(1/\spindex)$.  With this scaling, the
sufficient condition in Theorem~\ref{ThmSuff} reduces to $\numobs > C
\: \spindex \; \max \{ \log(\pdim - \spindex), \log
\frac{\pdim}{\spindex} \}$, whereas the necessary condition in
Theorem~\ref{ThmNec} reduces to $\numobs < C' \: \spindex \log
\frac{\pdim}{\spindex}$.  For many choices of sublinear sparsity
(e.g., $\spindex = \mathcal{O}(\sqrt{\pdim})$), we have $\log
\frac{\pdim}{\spindex} = \Omega(\log (\pdim - \spindex)) - o(1)$, so that we
can summarize the two conditions as a threshold of the order $\numobs
= \Theta(\spindex \log(\pdim - \spindex))$.  To compare with our
previous work~\cite{Wainwright06a_aller} on computationally tractable
methods, we established that $\ell_1$-constrained quadratic
programming (Lasso) has a threshold\footnote{Those
results~\cite{Wainwright06a_aller} allowed the minimum value to scale
as $\minval^2(\betastar) = f(\spindex)/\spindex$, where $f$ is any
function such that $\lim_{\spindex \rightarrow +\infty} f(\spindex) =
+\infty$.} for success/failure of order $\numobs = \Theta
\left(\spindex \log (\pdim - \spindex) \right)$, so that the Lasso
essentially achieves the information-theoretic bounds.
\myparagraph{Regime of linear sparsity} Next consider the regime of
linear sparsity, in which $\spindex = \alpha \pdim$ for some $\alpha
\in (0,1)$.  Considering first the sufficient conditions of
Theorem~\ref{ThmSuff}, we see that as long as $\minval^2(\betastar)
\spindex \rightarrow +\infty$, then $\numobs = \Theta(\pdim)$
observations are sufficient to ensure asymptotically reliable
recovery.  This information-theoretic condition should be compared
with our earlier analysis~\cite{Wainwright06a_aller} of
$\ell_1$-constrained quadratic programming (the Lasso); one
consequence of this work is that if \mbox{$\numobs < 2 \spindex \log
(\pdim -\spindex)$,} then the Lasso fails with probability converging
to one, even if $\minval^2(\betastar)$ stays bounded away from zero.
Given that $2 \spindex \log(\pdim-\spindex) \gg \Theta(\pdim)$ for
linear sparsity $\spindex = \alpha \pdim$, we see that there is a
substantial gap between the performance of the Lasso and the optimal
decoder in the linear sparsity regime. Thus, Theorem~\ref{ThmSuff}
raises the interesting question as to the existence of computationally
efficient techniques for asymptotically reliable recovery in the
regime of linear sparsity.

\section{Analysis}
\label{SecAnalysis}

This section is devoted to the proofs of Theorems~\ref{ThmSuff}
and~\ref{ThmNec}.  We begin by setting up some useful notation to be
used throughout the remainder of the paper.

\subsection{Notation and set-up}

For compactness in notation, let us use $\Amat$ to denote the $\numobs
\times \mdim$ matrix formed with the vectors $\arow_k = (\arow_{k1},
\arow_{k2}, \ldots, \arow_{k \mdim}) \in \real^\mdim$ as rows, and the
vectors $\acol_j = (\arow_{1j}, \arow_{2j}, \ldots, \arow_{\numobs
j})^T \in \real^\numobs$ as columns, as follows:
\begin{eqnarray}
\label{EqnAmatDefn}
\Amat & \defn & \begin{bmatrix} \arow_1^T \\ \arow_2^T \\ \vdots \\
\arow_\numobs^T \end{bmatrix} \; = \; \begin{bmatrix} \acol_1 & \acol_2
& \cdots & \acol_\mdim \end{bmatrix}.
\end{eqnarray}
Using $\Ysca$ and $\Wsca$ to denote the $\numobs$-dimensional
observation and noise vectors respectively, we can re-write our linear
observation model~\eqref{EqnLinearObs} in matrix-vector form as
follows:
\begin{eqnarray}
\label{EqnLinMatObs}
\Ysca & = & \Amat \betastar + \Wsca.
\end{eqnarray}
Given any subset $\Vset \subseteq \{1, \ldots, \pdim \}$, we use the
notation $\beta^*_\Vset$ to denote the $|\Vset|$-dimensional subvector
$\{ \beta^*_i, \; i \in \Vset \}$, and similarly for other vectors
(e.g., $\Ysca$, etc.).  In an analogous manner, we use $\Amatt{\Vset}$
to denote the $\numobs \times |\Vset|$ matrix with columns $\{
\acol_i, \; i \in \Vset \}$.  From herein, we assume without loss of
generality that $\sigma^2 = 1$, so that $\Wsca \sim N(0, I_{\numobs
\times \numobs}$) is simply a standard Gaussian vector.  (Note that
any scaling of $\sigma$ can be accounted for in the scaling of
$\betastar$, via the parameter $\minval(\betastar)$.)

In addition, we use the following standard notation for asymptotics of
real sequences $\{a_n\}$ and $\{b_n\}$: (i) $a_n = \mathcal{O}(b_n)$
means that $a_n \leq C b_n$ for some constant $C \in (0, \infty)$;
(ii) $a_n = \Omega(b_n)$ means that $a_n \geq C' b_n$ for some
constant $C' \in (0, \infty)$; (iii) $a_n = \Theta(b_n)$ is shorthand
for $a_n = \mathcal{O}(b_n)$ and $a_n = \Omega(b_n)$, and (iv) $a_n =
o(b_n)$ means that $a_n/b_n \rightarrow 0$.

\subsection{Proof of Theorem~\ref{ThmSuff}}
\label{SecSuff}

\myparagraph{Optimal decoding} We begin by describing the ``best''
decoder, that is optimal in terms of minimizing the probability of
error $\perr(\phi)$ over all decoding rules.  It is based on the
following real-valued function, defined on the subsets $\Uset \subset
\{1, \ldots, \mdim \}$, as
\begin{equation}
\OptFun(\Uset; \Ysca, \Amat, \betastar) \; = \; \arg
\min_{\Beta{\Uset}} \left \{\|Y - \Amatt{\Uset} \Beta{\Uset} \|_2^2
\right \}.
\end{equation}
We frequently write $\OptFun(\Uset)$ as a shorthand; note that this
value corresponds to the error associated with the best estimator of
$Y$ that lies in $\Range(\Amatt{\Uset})$.  The optimal decoder chooses
the best subset $\Sbest$ based on the minimal value of this error,
ranging over all subsets $\Uset$ of size $\spindex$:
\begin{eqnarray}
\label{EqnOptDecoder}
\Sbest \; = \; \decopt(Y) & \defn & \arg \min_{|\Uset | = \spindex }
\OptFun(\Uset; \Ysca, \Amat, \betastar).
\end{eqnarray}
Note that by symmetry, the error probability \mbox{$\Prob[\Sbest \neq
\Sset \, \mid \, \Sset]$} is in fact the same regardless of which
underlying set $\Sset$ acts as the true one.  Consequently, we can
view the choice of $\Sset$ as fixed (and hence non-random), and write
\begin{equation}
\perr(\phi) \; = \; \Prob[\phi(Y) \neq \Sset],
\end{equation}
which should now be understood as an unconditional probability (with
$\Sset$ fixed).

\myparagraph{Analysis of error probability} Consider the difference
$\Metric(\Uset) \defn \OptFun(\Uset) - \OptFun(\Sset)$ between the
reconstruction error $\OptFun(\Sset)$ using the true subset $\Sset$,
versus the error $\OptFun(\Uset)$ candidate subset $\Uset$.  For any
subset $\Uset$ such that $\Amatt{\Uset}$ is full rank, define the
$\numobs \times \numobs$ matrices
\begin{subequations}
\begin{eqnarray}
\Projmat{\Uset} & \defn & \Gfullmat{\Amatt{\Uset}}, \qquad \qquad
\mbox{and} \\
\Projmatc{\Uset} & \defn & I_{\numobs \times \numobs} -
\Gfullmat{\Amatt{\Uset}}.
\end{eqnarray}
\end{subequations}
Note that $\Projmat{\Uset}$ and $\Projmatc{\Uset}$ are both orthogonal
projection matrices, associated with the $\spindex$-dimensional range
space $\Range(\Amatt{\Uset})$ and $(\numobs-\spindex)$-dimensional
nullspace $\Ker(\Amatt{\Uset})$ respectively.  With these definitions,
we state the following result (see Appendix~\ref{AppAlgebra} for a
proof):
\blems
\label{LemAnalyzeForm}
For a given vector $\betastar$ with support $\Sset$, the optimal
decoder declares $\Uset$ over $\Sset$ if and only if the random
variable
\begin{eqnarray}
\label{EqnAnalyzeForm}
\Metric(\Uset) & = & \left \| \Projmatc{\Uset} \left(\Amatt{\SbackU}
\betastar_{\SbackU} + \Wsca \right) \right \|^2 - \left \|
\Projmatc{\Sset} \Wsca \right \|^2.
\end{eqnarray}
is negative.
\elems
\noindent Overall, the optimal decoder fails if and only if at least
one $\Uset$ (with cardinality $|\Uset| = \spindex$) is preferable to
$\Sset$; consequently, the probability of error can be written as
\begin{eqnarray}
\label{EqnErrorProbOne}
\mprob[\Sbest \neq \Sset] & = & \mprob \Big[\bigcup_{\Uset \neq \Sset,
    \; |\Uset| = \spindex} \{ \Metric(\Uset) < 0 \} \Big].
\end{eqnarray}
In order to analyze this error probability, we begin by considering
the range of possible integers $k \defn |\Sset \backslash \Uset|$,
corresponding to the complement of the overlap.  The following lemma
characterizes the exponential decay rates of the random variable
$\Metric(\Uset)$:
\blems
\label{LemExpUpper}
For fixed $k$ (with $1 \leq k \leq \spindex$), we have for any $\Uset$
with $|\SbackU| = k$,
\begin{eqnarray}
\label{EqnExpUpper}
\prob[\Metric(\Uset) < 0 ] & \leq & \exp \left \{\frac{-(\numobs -
  \spindex)\|\betastar_{\SbackU}\|^2}{12
  \left(\|\betastar_{\SbackU}\|^2 + 4\right)} \right \} + 2 \, \exp
  \left \{ -\frac{k}{4} \left[-1 + \frac{1}{4}(\numobs - \spindex)
  \frac{ \|\betastar_{\SbackU}\|^2}{k} \right]^2 \right \}.
\end{eqnarray}
\elems
\spro
We begin by conditioning on the Gaussian noise vector $\Wsca$. Since
each element of $\Amatt{\SbackU}$ is standard normal, each entry of
the random vector $\Amatt{\SbackU} \betastar_{\SbackU}$ is zero-mean
Gaussian with variance $\|\betastar_{\SbackU}\|^2$.  Consequently, if
we rescale by the standard deviation, then the random vector
\begin{equation*}
\|\betastar_{\SbackU}\|^{-1} \; \left(\Amatt{\SbackU}
\betastar_{\SbackU} + \Wsca\right)
\end{equation*}
is an $\numobs$-dimensional Gaussian random vector with independent
entries, each with with unit variance, and mean vector $\Wsca$.
Applying the orthogonal transform $\Projmatc{\Uset}$ reduces the
number of degrees of freedom to $(\numobs - \spindex)$, so that we
conclude that
\begin{equation*}
\|\betastar_{\SbackU}\|^{-2} \; \left\| \Projmatc{\Uset} \left(
\Amatt{\SbackU} \betastar_{\SbackU} + \Wsca \right) \right \|^2
\end{equation*}
is a non-central $\chi^2$ variate with $d = \numobs - \spindex$
degrees of freedom, and non-centrality parameter $\nu =
\|\betastar_{\SbackU} \|^{-2} \|\Projmatc{\Uset} \Wsca\|^2$.  With
these choices of $(d, \nu)$, we have
\begin{eqnarray*}
\mprob[\Metric(\Uset) < 0 \, \mid \, \Wsca] & = & \mprob
\left[\chi^2(d, \nu) < t \right]
\end{eqnarray*}
where we have set $t \defn \frac{\|\Projmatc{\Sset}
\Wsca\|^2}{\|\betastar_{\SbackU} \|^2}$ for shorthand.  Thus,
conditioned on $\Wsca$, our problem reduces to bounding the tail of a
non-central $\chi^2$ variate.  In Appendix~\ref{AppChiTail}, we state
some known tail bounds~\cite{Birge01} on such variates, which we use
here.  In order to apply these bounds, we condition on the following
``good event'', defined in terms of $\Wsca$
\begin{eqnarray*}
\goodW & = & \left\{ \left | \frac{ \|\Projmatc{\Uset} \Wsca\|^2 -
  \|\Projmatc{\Sset} \Wsca \|^2}{\|\betastar_{\SbackU}\|^2} \right |
  \leq \frac{\numobs- \spindex}{2} \right \} \bigcap \left \{
  \|\Projmatc{\Uset} \Wsca\|^2 \leq 2 (\numobs - \spindex) \right \}.
\end{eqnarray*}
Note that the first event defining $\goodW$ ensures that
\begin{eqnarray}
\label{EqnKeyCond}
d + \nu - t & = & (\numobs - \spindex) +
\frac{1}{\|\betastar_{\SbackU}\|^2} \left(\|\Projmatc{\Uset} \Wsca
\|^2 - \|\Projmatc{\Sset} \Wsca \|^2 \right) \; \geq \;
\frac{\numobs-\spindex}{2} \;\geq \; 0.
\end{eqnarray}
Consequently, conditioned on $\goodW$, we may set $x \defn
\frac{\left(d + \nu - t\right)^2}{4 (d + 2\nu)}$ in
equation~\eqref{EqnNoncentB} to obtain the upper bound
\begin{eqnarray}
\log \mprob[ \Metric(\Uset) < 0 \; \mid \; \goodW] & \leq &
-\frac{\left (d + \nu - t \right)^2}{4 (d + 2\nu)} \nonumber \\
& = & -\frac{ \left ( [\numobs - \spindex] + \frac{\|\Projmatc{\Uset}
\Wsca\|^2 - \|\Projmatc{\Sset} \Wsca\|^2}{\|\betastar_{\SbackU}\|^2}
\right)^2}{4 \left( [\numobs-\spindex] + 2 \frac{\|\Projmatc{\Uset}
\Wsca\|^2}{\|\betastar_{\SbackU} \|^{2}} \right)} \nonumber \\
& = & - \left(\numobs - \spindex \right) \; \frac{ \left ( 1 +
\frac{\|\Projmatc{\Uset} \Wsca\|^2 - \|\Projmatc{\Sset}
\Wsca\|^2}{(\numobs-\spindex) \; \|\betastar_{\SbackU}\|^2}
\right)^2}{4 \left(1+ 2 \frac{\|\Projmatc{\Uset} \Wsca\|^2}{(\numobs -
\spindex) \; \|\betastar_{\SbackU} \|^{2}} \right)} \nonumber \\
& \stackrel{(b)}{\leq} & - \left(\numobs - \spindex \right) \;
\frac{1/2}{4\left(1+ 4/\|\betastar_{\SbackU}\|^2 \right)} \nonumber \\
\label{EqnTermOne}
& = & - \left(\numobs - \spindex \right) \;
\frac{\|\betastar_{\SbackU}\|^2}{8 \left(\|\betastar_{\SbackU}\|^2 +
4\right)},
\end{eqnarray}
where inequality (b) makes use of the second event defining $\goodW$.

We complete the proof by observing that
\begin{eqnarray}
\label{EqnOverall}
\mprob[ \Metric(\Uset) < 0] & \leq & \mprob[ \Metric(\Uset) < 0 \;
\mid \; \goodW] + \mprob[\goodW^c],
\end{eqnarray}
so that it suffices to upper bound $\mprob[\goodW^c]$.  By union
bound, we have
\begin{eqnarray}
\label{EqnUnion}
\mprob[\goodW^c] & \leq & \mprob \left[\left | \frac{
  \|\Projmatc{\Uset} \Wsca\|^2 - \|\Projmatc{\Sset} \Wsca
  \|^2}{\|\betastar_{\SbackU}\|^2} \right | \geq \frac{\numobs-
  \spindex}{2} \right] + \mprob \left[\|\Projmatc{\Uset} \Wsca\|^2
  \geq 2 (\numobs - \spindex) \right].
\end{eqnarray}
Since $\|\Projmatc{\Uset} \Wsca\|^2$ is a central $\chi^2$ with
$(\numobs - \spindex)$ degrees of freedom, we may apply the tail
bounds from Appendix~\ref{AppChiTail} to conclude that 
\begin{eqnarray}
\label{EqnTermThree}
\mprob\left[\|\Projmatc{\Uset} \Wsca\|^2 \geq 2 (\numobs - \spindex)
\right] & \leq & \exp(-(\numobs - \spindex)/12).
\end{eqnarray}
Turning to the first term on the RHS on equation~\eqref{EqnUnion}, we
observe that
\begin{eqnarray*}
 \|\Projmatc{\Uset} \Wsca\|^2 - \|\Projmatc{\Sset} \Wsca \|^2 & = &
  \|\Projmat{\Uset} \Wsca\|^2 - \|\Projmat{\Sset} \Wsca \|^2 \; \edist
  \; \sum_{i \in \UbackS} Z_i^2 - \sum_{j \in \SbackU} Z_j^2,
\end{eqnarray*}
where $\{Z_i, Z_j \}$ are i.i.d. standard normal variates.  Now if the
difference $\sum_{i \in \UbackS} Z_i^2 - \sum_{j \in \SbackU} Z_j^2$
is to exceed $\frac{1}{2}(\numobs - \spindex)
\|\betastar_{\SbackU}\|^2$, then at least one of the terms must exceed
$\frac{1}{4}(\numobs - \spindex) \|\betastar_{\SbackU}\|^2$.
Moreover, we observe that $\sum_{j \in \SbackU} Z_j^2$ is $\chi^2_k$,
where $k = |\SbackU|$.  Hence, we have
\begin{eqnarray*}
\log \mprob \left[\left | \frac{ \|\Projmatc{\Uset} \Wsca\|^2 -
  \|\Projmatc{\Sset} \Wsca \|^2}{\|\betastar_{\SbackU}\|^2} \right |
  \geq \frac{\numobs- \spindex}{2} \right] & \leq & \log 2 \, \mprob
  \left[ \frac{\chi^2_k}{k} \geq \frac{1}{4}(\numobs - \spindex)
  \frac{\|\betastar_{\SbackU}\|^2}{k} \right] \\
& = & \log 2 \mprob \left[\chi^2_k -k \geq k \left\{ -1 +
  \frac{1}{4}(\numobs - \spindex) \frac{\|\betastar_{\SbackU}\|^2}{k}
  \right \} \right] \\
& \leq & -\frac{k}{4} \left[-1 + \frac{1}{4}(\numobs - \spindex)
 \frac{ \|\betastar_{\SbackU}\|^2}{k} \right]^2 + \log 2,
\end{eqnarray*}
where we have used the upper bound~\eqref{EqnCleanUpCent} from
Appendix~\ref{AppChiTail} with $x \defn \frac{k}{4} \left( -1 +
\frac{1}{4}(\numobs - \spindex) \frac{\|\betastar_{\SbackU}\|^2}{k}
\right)^2$ in the final inequality.

\fpro

\myparagraph{Weakened but simpler bound} In order to make further
progress, we simplify the bound~\eqref{EqnExpUpper} from
Lemma~\ref{LemExpUpper}, at the expense of weakening it, by noting
that for all $k \geq 1$, we have $\|\betastar_{\SbackU}\|^2 \geq k \,
\minval^2(\betastar)$, so that
\begin{eqnarray}
\label{EqnWeakOne}
\prob[\Metric(\Uset) \leq 0 ] & \leq & \exp \left \{\frac{-(\numobs -
  \spindex) k \ \minval^2(\betastar)}{12 \left(k \,
  \minval^2(\betastar) + 4\right)} \right \} + 2 \, \exp \left
  \{-\frac{k}{4} \left[\frac{\numobs - \spindex}{4} \:
  \minval^2(\betastar) - 1 \right]^2 \right \}.
\end{eqnarray}
The advantage of this weakened bound is that it is independent of the
subset $\Uset$, and depends only on the parameter $k = |\SbackU|$.

From this weakened bound~\eqref{EqnWeakOne}, we see the necessity (at
least for this analysis) of the requirement $(\numobs - \spindex) \,
\minval^2(\betastar) \rightarrow +\infty$, so that the second error
term decays asymptotically.  Under this requirement, we have (for
sufficiently large $\numobs$) that the second error exponent can be
bounded as
\begin{eqnarray*}
-\frac{k}{4} \left[\frac{\numobs - \spindex}{4} \:
  \minval^2(\betastar) - 1 \right]^2 & \leq & -\frac{k}{12} \left[
  \frac{\numobs - \spindex}{4} \: \minval^2(\betastar) - 1\right] \\
 & \leq & -\frac{k}{4} \frac{\numobs - \spindex}{8} \:
 \minval^2(\betastar) \\
& \leq & \frac{-(\numobs - \spindex) k \, \minval^2(\betastar)}{12
 \left(k \, \minval^2(\betastar) + 8\right)}.
\end{eqnarray*}
The first error exponent is also upper bounded by this same quantity,
so that we can simplify the upper bound to
\begin{eqnarray}
\label{EqnWeakTwo}
\prob[\Metric(\Uset) \leq 0 ] & \leq & 3\, \exp \left \{
  \frac{-(\numobs - \spindex) k \, \minval^2(\betastar)}{12 \left(k \,
  \minval^2(\betastar) + 8\right)} \right \}.
\end{eqnarray}
Denote by $N(k)$ the number of subsets $\Uset$ of size $\spindex$,
with overlap exactly equal to $k$.  A standard counting argument
yields that, for each $k$ with $1 \leq k \leq \spindex$, there are
\begin{eqnarray}
\label{EqnExactNK}
N(k) & = & {\spindex \choose k} \, {\pdim - \spindex \choose k}
\end{eqnarray}
such subsets.  Using this simple bound~\eqref{EqnWeakTwo} and union
bound applied to the representation~\eqref{EqnErrorProbOne}, we can
upper bound the error probability as
\begin{eqnarray}
\label{EqnFinalUnion}
\mprob[\Sbest \neq \Sset] & \leq & 3 \sum_{k=1}^\spindex {\spindex
\choose k} \, {\pdim - \spindex \choose k} \; \exp \left \{
\frac{-(\numobs - \spindex) k \, \minval^2(\betastar)}{12 \left(k \,
\minval^2(\betastar) + 8\right)} \right \}.
\end{eqnarray}

\myparagraph{Analysis of the upper bound} We now analyze the upper
bound~\eqref{EqnFinalUnion}; in particular, our goal is to derive
sufficient conditions for each of the terms in the summation to vanish
asymptotically.  In order to deal with the binomial coefficients, we
make use of the bounds (see Appendix~\ref{AppBinBound})
\begin{equation}
\log {\spindex \choose k} \leq k \log \frac{\spindex e}{k}, \qquad
\mbox{and} \qquad
\log {\pdim - \spindex \choose k} \leq k \log \frac{(\pdim-\spindex)
e}{k}.
\end{equation}
Applying these two bounds, we conclude that the (logarithm of the)
$k^{th}$ term is upper bounded by
\begin{equation*}
k \left [ 2 + \log \frac{\spindex}{k} + \log \frac{\pdim -
\spindex}{k} \right] - \frac{(\numobs - \spindex) k \,
\minval^2(\betastar)}{12 \left(k \, \minval^2(\betastar) + 8\right)}.
\end{equation*}
Requiring this term to be negative asymptotically is equivalent to
having
\begin{eqnarray}
(\numobs - \spindex) & \geq & \frac{12 \left(k \minval^2(\betastar) +
8\right)}{k \, \minval^2(\betastar)} \; k \left [ 2 + \log
\frac{\spindex}{k} + \log \frac{\pdim - \spindex}{k} \right] \nonumber
\\
\label{EqnInterBound}
& = & 12 \left(k \, + \frac{8}{\minval^2(\betastar)} \right) \left \{
2 + \log \frac{\spindex}{k} + \log \frac{\pdim - \spindex}{k}
\right\}.
\end{eqnarray}

In order to understand the behavior of this lower bound, we consider
$k$ in two distinct regimes: 
\bit
\item On one hand, if $k = \mycfrac \spindex$ for some $\mycfrac \in
(0,1)$, then the second term on the RHS of the
bound~\eqref{EqnInterBound} is dominated by the term $\log \frac{\pdim
- \spindex}{\mycfrac \spindex} = \Omega(\log \frac{\pdim}{\spindex})$,
so that the overall lower bound is dominated by $\max \{\spindex,
\minval^{-2}(\betastar)\} \log(\pdim/\spindex)$.
\item On the other hand, if $k = o(\spindex)$, the lower bound is
dominated by the maximum of linear growth $\spindex$, and the quantity
\mbox{$\minval^{-2}(\betastar) \, \log(\pdim - \spindex)$.}
\eit
Overall, we conclude that the condition
\begin{eqnarray}
\label{EqnFinalSublinCond}
\numobs & > & C \; \max\left \{ \spindex \log (\pdim/\spindex), \;
\frac{1}{\minval^2(\betastar)} \log(\pdim - \spindex) \right \},
\end{eqnarray}
for some constant $C > 0$ is sufficient in order to achieve
asymptotically reliable recovery, as claimed in Theorem~\ref{ThmSuff}.

\subsection{Proof of Theorem~\ref{ThmNec}}
\label{SecNec}

We now turn to the proof of the necessary conditions given in
Theorem~\ref{ThmNec}.  

\myparagraph{Fano method} Our analysis is based on a well-known lower
bound on the probability of error in a multiway hypothesis testing
problem in terms of Kullback-Leibler divergences. In the
non-parametric statistics
literature~\cite{Hasminskii78,IbrHas81,Yu97}, this approach is
referred to as the Fano method, since the bound is a corollary of
Fano's inequality from information theory~\cite{Cover}.  Here we state
and make use of the following variant~\cite{Yatracos88}:
\blems
Consider a family of $\Numset$ distributions $\{\prob_1, \ldots,
\prob_\Numset \}$.  Then the average probability of error in
performing in a hypothesis test over this family is lower bounded as
\begin{eqnarray*}
\perr & \geq & 1 - \frac{ \frac{1}{\Numset^2} \sum \limits_{i,j
=1}^\Numset D(\prob_i \, \| \, \prob_j) + \log 2}{\log \left(\Numset-1
\right)},
\end{eqnarray*}
where $D(\prob_i \, \| \, \prob_j)$ denotes the Kullback-Leibler
divergence between distributions $\prob_i$ and $\prob_j$.
\elems

\myparagraph{Restricted problem} Consider the collection of all
$\Numset = {\pdim \choose \spindex}$ subsets of size $\spindex$ chosen
from $\{1, \ldots, \pdim\}$.  In order to produce lower bounds, we
analyze the behavior of the optimal decoder for a restricted problem,
in which we assume that for any fixed support $\Sset$, it is known
\emph{a priori} that $\betastar_i = \minval(\betastar)$ for all
indices $i \in \Sset$.  (Recall that $\minval(\betastar)$ is the
minimum absolute value of entries in the support of $\betastar$.)
This problem is simply an $\Numset$-way hypothesis testing problem, in
which the observation under the hypothesis associated with subset
$\Uset$ takes the form 
\begin{eqnarray}
\label{EqnModObs}
\Ysca & = &  \Amatt{\Uset} \mybetavec + \Wsca,
\end{eqnarray}
where $\mybetavec = \minval(\betastar) \myones_\spindex$ is a rescaled
$\spindex$-vector of ones, and $\Wsca \sim N(0, I_{\numobs \times
\numobs})$.

Let us index the collection of all $\spindex$-sized subsets with $i =
1, 2, \ldots, \Numset$, and use $\Uset[i]$ to denote the corresponding
support.  For each index $i$, let $\prob_i$ denote the multivariate
Gaussian distribution with mean $\Amatt{\Uset[i]} \mybetavec$ and
covariance matrix $I_{\numobs \times \numobs}$; note that $\prob_i$ is
simply the class-conditional distribution of $\Ysca$ under the
hypothesis $\Uset[i]$.  Moreover, the Kullback-Leibler divergence
between any such pair is given by $D(\prob_i \, \| \, \prob_j) =
\frac{1}{2} \|\Amatt{\Uset[i]} \mybetavec - \Amatt{\Uset[j]}
\mybetavec \|_2^2$, so that the corresponding Fano bound takes the
form
\begin{eqnarray*}
\perr & \geq & 1 - \frac{1}{2} \frac{ \frac{1}{ N^2} \sum_{i, j
=1}^{N} \|\Amatt{\Uset[i]} \mybetavec - \Amatt{\Uset[j]} \mybetavec
\|_2^2 + 2 \log 2}{\log [N-1]}.
\end{eqnarray*}

\myparagraph{Upper bounds via concentration} Thus, in order to ensure
that $p_e$ stays bounded away from zero, we need to (upper) bound the
quantity $\frac{1}{2} \frac{1}{\Numset^2} \sum_{i,
j=1}^{\Numset} \|\Amatt{\Uset[i]} \mybetavec -
\Amatt{\Uset[j]} \mybetavec \|_2^2 \big / \log [\Numset-1]$
away from one.  For a given pair of subsets $(\Uset, \Vset)$ in our
collection, consider the random variable \mbox{$Z_{\Uset, \Vset} \defn
\|\Amatt{\Uset} \mybetavec - \Amatt{\Vset} \mybetavec \|_2^2$.}  A
little calculation shows that $Z_{\Uset, \Vset} \sim \gamvar(\Uset,
\Vset) \chi^2_\numobs$, where
\begin{equation}
\label{EqnDefnGamvar}
\gamvar(\Uset, \Vset) = 2 \, \minval^2(\betastar) \, \left (\spindex -
|\Uset \cap \Vset| \right).
\end{equation}
The following result bounds the upper tail behavior of the random
variable $Z = \frac{1}{\Numset^2} \sum_{\Uset \neq \Vset} Z_{\Uset,
\Vset}$.
\blems
\label{LemConcen}
The tail of $Z$ obeys the bound
\begin{eqnarray*}
\prob \left[Z \geq 4 \minval^2(\betastar) \spindex \numobs \right] &
\leq & \frac{1}{2}.
\end{eqnarray*}
\elems
\noindent Using this lemma (see Appendix~\ref{AppConcen} for a proof
of this claim), we are guaranteed that at least $1/2$ of the Gaussian
ensembles satisfy the upper bound
\begin{eqnarray}
\label{EqnKeyQuant}
\frac{1}{2} \; \frac{\frac{1}{\Numset^2} \sum \limits_{i,j
=1}^{\Numset} D(\prob_i \, \| \, \prob_j)}{\log
[\Numset-1]} & = & \frac{1}{2} \;
\frac{\frac{1}{\Numset^2} \sum_{\Uset \neq \Vset} Z_{\Uset,
\Vset}}{\log[\Numset-1]} \; \leq \;
 \frac{4 \minval^2(\betastar) \spindex \numobs}{\log
[\Numset-1]}.
\end{eqnarray}
Hence, as long as the quantity~\eqref{EqnKeyQuant} remains bounded
from above away from one, the Fano bound implies that the probability
of error averaged over the whole ensemble will remain bounded away
from zero.  Consequently, we obtain the necessary condition that
\begin{eqnarray*}
\numobs & > & \frac{\log [\Numset-1]}{4 \minval^2(\betastar)
\spindex}
\end{eqnarray*}
for reliable recovery with probability one asymptotically. To obtain a
more transparent bound, we first lower bound $N$ via $\log
[N-1] \geq \frac{1}{2} \log \Numset$, and then
further via
\begin{eqnarray*}
\frac{1}{2} \log \Numset & = & \frac{1}{2} \log {\pdim \choose
\spindex} \; \geq \; \frac{1}{2} \spindex \; \log
\frac{\pdim}{\spindex},
\end{eqnarray*}
as stated in Appendix~\ref{AppBinBound}.  Consequently, we obtain the
necessary condition
\begin{eqnarray}
\numobs > \Omega \left(\frac{1}{\spindex \; \minval^2(\betastar)}
\spindex \log \frac{\pdim}{\spindex} \right),
\end{eqnarray}
as stated in Theorem~\ref{ThmNec}.

\section{Conclusion}
\label{SecDiscussion}

In this paper, we have analyzed the information-theoretic limits of
the sparsity recovery problem for the linear observation
model~\eqref{EqnLinearObs} with measurement vectors drawn from the
standard Gaussian ensemble.  We have established both lower and upper
bounds on the number of observations $\numobs$ as a function of the
model dimension $\pdim$ and sparsity index $\spindex$ that are
required for asymptotically reliable recovery.

There are a variety of open questions raised by our analysis. First,
while our upper and lower bounds are essentially matching for certain
regimes of scaling (e.g., sublinear sparsity with the minimum
$\minval^2(\betastar) = \Theta(1/\spindex)$), it is likely that the
analysis can be tightened in other regimes.  In particular, the
analysis of the necessary conditions (see proof of
Theorem~\ref{ThmNec}) involves some slack since it is based on
analyzing a very restricted ensemble.  Second, our results (in
particular, a corollary of Theorem~\ref{ThmSuff}) reveal that with the
sparsity index scaling linearly ($\spindex = \alpha \pdim$ for some
$\alpha \in (0,1)$), as long, as the minimum value
$\minval^2(\betastar)$ decays sufficiently slowly, then asymptotically
reliable recovery is possible with only a linear number of
observations (i.e., $\numobs = \beta \pdim$ for some $\beta > 0$).
Since our previous work~\cite{Wainwright06a_aller} established that
the Lasso ($\ell_1$-constrained quadratic programming) cannot achieve
reliable recovery in this particular $(\numobs, \pdim, \spindex)$
regime, it remains to determine a computationally tractable method
that approaches such performance in the regime of linear sparsity.
Third, whereas the current analysis has focused on a very special
class of Gaussian ensemble, the analysis given here could be extended
to a broader class of measurement ensembles.

\subsection*{Acknowledgements}  This work was partially supported by
NSF CAREER Award CCF-0545862, NSF Grant DMS-0605165, and an Alfred P.
Sloan Foundation Fellowship.  We thank Peter Bickel for helpful
discussions and pointers.


\appendix

\section{Proof of Lemma~\ref{LemAnalyzeForm}}
\label{AppAlgebra}

We begin by showing that for any subset $\Uset$ for which
$\Amatt{\Uset}$ is full rank, the function $\OptFun$ has the
equivalent form $\OptFun(\Uset) = \|\Projmatc{\Uset} \Ysca \|_2^2$.
Under the given rank condition, the linear least squares estimator of
$\betastar_\Uset$ is given by $\estim{\beta}_{\Uset} =
\MatInv{\Amatt{\Uset}} \Amatt{\Uset}^T \Ysca$.  Noting that
$\Amatt{\Uset} \estim{\beta}_{\Uset} = \Projmat{\Uset} \Ysca$, we
substitute into the quadratic norm and expand, thereby obtaining
\begin{eqnarray*}
\OptFun(\Uset) & = & \| \Ysca - \Amatt{\Uset} \estim{\Beta{\Uset}}
\|_2^2 \; = \; \|(I - \Projmat{\Uset}) \Ysca\|_2^2 \; = \;
\|\Projmatc{\Uset} \Ysca \|_2^2
\end{eqnarray*}
as claimed.  Lastly, to establish equation~\eqref{EqnAnalyzeForm}, we
note that
\begin{eqnarray*}
\OptFun(\Uset) & = & \|\Projmatc{\Uset} \left( \Amatt{\Sset}
\betastar_\Sset + \Esca \right) \|_2^2 \; = \; \|\Projmatc{\Uset}
\left( \Amatt{\SbackU} \betastar_\SbackU + \Esca \right) \|_2^2,
\end{eqnarray*}
since $\Projmatc{\Uset} v = 0$ for any vector $v$ belonging to the
range of $\Amatt{\Uset}$.



\section{Proof of Lemma~\ref{LemConcen}}
\label{AppConcen}
Note that $Z = \frac{1}{\Numset^2} \sum_{\Uset, \Vset} Z_{\Uset,
\Vset}$ is a rescaled sum of a total number $\Numset^2$ variables
(neither independent nor identically distributed).  However, since $Z$
is a non-negative random variable, we may apply Markov's inequality
for any $t > 0$ to conclude that
\begin{eqnarray}
\label{EqnMarkov}
\prob\left[ Z \geq t \right] & \leq & \frac{\Exs[Z]}{t}.
\end{eqnarray}
Since each $Z_{\Uset, \Vset}$ has distribution $\gamvar(\Uset, \Vset)
\chi^2_\numobs$, we have $\Exs[Z_{\Uset, \Vset}] = \gamvar(\Uset,
\Vset) \numobs$.  From equation~\eqref{EqnDefnGamvar}, we note that
$\gamvar(\Uset, \Vset) \leq 2 \minval^2(\betastar) \spindex$, and hence
\begin{eqnarray*}
\Exs[Z] & \leq & \max_{\Uset \neq \Vset} \left(\gamvar(\Uset,
\Vset)\right) \numobs \; = \; 2 \minval^2(\betastar) \spindex \numobs,
\end{eqnarray*}
Hence setting $t = 4 \minval^2(\betastar) \spindex \numobs$ in the
bound~\eqref{EqnMarkov} yields the claim.


\section{Bounds on binomial coefficients}
\label{AppBinBound}

Although more refined results are certainly possible, we make frequent
use of the following crude bounds on the binomial coefficients
\begin{equation}
\label{EqnBinBound}
\left(\frac{n}{k} \right)^k \; \leq \; {n \choose k} \; \leq \;
\left(\frac{n \, e}{k} \right)^k.
\end{equation}

\section{Tail bounds for chi-square variables}
\label{AppChiTail}
The following large-deviations bounds for centralized $\chi^2$ are
taken from Laurent and Massart~\cite{LauMas98}.  Given a centralized
$\chi^2$-variate $X$ with $d$ degrees of freedom, then for all $x \geq
0$,
\begin{subequations}
\begin{eqnarray}
\label{EqnCleanUpCent}
\mprob \left[X - d \geq 2 \sqrt{d x} + 2x \right] \; \leq \; \mprob
\left[X - d \geq 2\sqrt{d x} \right] & \leq & \exp(-x), \qquad \mbox{and} \\
\label{EqnCleanDownCent}
\mprob \left[X - d \leq -2 \sqrt{d x} \right] & \leq & \exp(-x).
\end{eqnarray}
\end{subequations}
More generally, the analogous tail bounds for \emph{non-central}
$\chi^2$, taken from Birg\'{e}~\cite{Birge01}, can be established via
the Chernoff bound.  Let $X$ be a non-central $\chi^2$ variable with
$\df$ degrees of freedom and non-centrality parameter $\noncent \geq
0$.  Then for all $x > 0$,
\begin{subequations}
\begin{eqnarray}
\label{EqnNoncentA}
\prob \left[ X \geq (\df + \noncent) + 2 \sqrt{(\df + 2 \noncent) x} +
2 x \right] & \leq & \exp(-x), \qquad \mbox{and} \\
\label{EqnNoncentB}
\prob \left[ X \leq (\df + \noncent) -2 \sqrt{(\df + 2 \noncent) x}
\right] & \leq  & \exp(-x).
\end{eqnarray}
\end{subequations}


\end{document}